\documentclass[11pt]{amsart}
\usepackage{graphicx}
\usepackage{graphics}
\usepackage{amsmath,overpic}
\usepackage{amscd}
\usepackage{latexsym}
\usepackage[all]{xy}
\begin{document}
\textwidth 5.5in
\textheight 8.3in
\evensidemargin .75in
\oddsidemargin.75in

\newtheorem{lem}{Lemma}[section]
\newtheorem{fact}[lem]{Fact}
\newtheorem{conj}[lem]{Conjecture}
\newtheorem{defn}[lem]{Definition}
\newtheorem{thm}[lem]{Theorem}
\newtheorem{cor}[lem]{Corollary}
\newtheorem{prob}[lem]{Problem}
\newtheorem{claim}[lem]{Claim}
\newtheorem{main}{Main theorem}
\newtheorem{exm}[lem]{Example}
\newtheorem{rmk}[lem]{Remark}
\newtheorem{que}[lem]{Question}
\newtheorem{prop}[lem]{Proposition}
\newtheorem*{observation}{Observation}
\newtheorem{clm}[lem]{Claim}
\newcommand{\p}[3]{OS_{p,#1}^{#2}(#3)}
\def\Z{\mathbb Z}
\def\ker{\text{Ker}}
\def\F{\mathbb F}
\def\C{\mathcal{C}}
\def\D{\mathcal{D}}
\def\P{\mathcal{P}}
\def\R{\mathbb R}
\def\bN{\mathbb N}
\def\X{\mathcal{X}}
\def\g{\overline{g}}
\def\bZ{\mathbb Z}
\def\bF{\mathbb F}
\def\odots{\reflectbox{\text{$\ddots$}}}
\newcommand{\tg}{\overline{g}}
\def\proof{{\bf Proof. }}
\def\ee{\epsilon_1'}
\def\ef{\epsilon_2'}
\title{Non-existence theorems on infinite order corks}
\author{Motoo Tange}
\thanks{The author was partially supported by JSPS KAKENHI Grant Number 26800031}
\address{Institute of Mathematics, University of Tsukuba,
 1-1-1 Tennodai, Tsukuba, Ibaraki 305-8571, Japan}
\email{tange@math.tsukuba.ac.jp}
\date{\today}
\maketitle
\begin{abstract}
Suppose that $X,X'$ are simply-connected closed exotic 4-manifolds.
It is well-known that $X'$ is obtained by an order 2 cork twist of $X$.
We give an infinite exotic family of 4-manifolds not generated by any infinite order cork.
This is the first example admitting such a condition.
We prove a necessary condition of 4-dimensional OS-invariants for a family to be generated by an infinite order cork
and give non-contractible relatively exotic 4-manifolds that are never induced by any cork.
Furthermore, we prove an estimate of the number of OS-invariants of 4-manifolds generated by a cork.
\end{abstract}

\section{Introduction}
\label{intro}
\subsection{Motivation and some results}
\label{Corktwist}
We assume that all manifolds in this paper are smooth and oriented.
`Exotic' means that two manifolds are homeomorphic but non-diffeomorphic each other.
For two exotic 4-manifolds, the following fact is well-known.
\begin{fact}[\cite{Mat},\cite{CFHS}]
\label{Ma}
Let $X$ and $X'$ be two simply-connected, closed exotic smooth 4-manifolds.
Then, there exists a contractible 4-manifold $C$ embedded in $X$ by an embedding $e$ such that $X'=(X-e(C))\cup_\tau C$ and $\tau^2=\text{id}$.
\end{fact}
Here $\tau^n$ presents the $n$-fold composition $\tau\circ \cdots \circ \tau$ of a map $\tau$.
The notation $\cup_\tau$ stands for gluing two manifolds via a diffeomorphism $\tau$ on the boundaries.
Through this paper, we denote the resulting manifold $(X-e(C))\cup_{\tau} C$ by $X(C,e,\tau)$.
We omit $e$ from the notation if we can understand the embedding from the context.

As mentioned in \cite{AM}, we may assume that the contractible 4-manifold $C$ in Fact~\ref{Ma} is a Stein manifold.
The contractible 4-manifold satisfying Fact~\ref{Ma} is called a cork, and the general definition is in Section~\ref{copl}.
We call the fact a {\it cork existence theorem}.
Fact~\ref{Ma} means that any two exotic, closed, simply-connected 4-manifolds
$\mathcal{X}=\{X_0,X_1\}$ are related by a surgery $X_i\leadsto X_i(C,e,\tau)=X_{i+1}$, where the index set of $i$ is ${\mathbb Z}/2{\mathbb Z}$.

Our main questions are the following:
\begin{itemize}
\item For any general family $\mathcal{X}$ of closed exotic 4-manifolds, 
can we find a triple $(C,e,\tau)$ satisfying $\mathcal{X}\subset \{X(C,e,\tau^i)|i\in{\mathbb Z}\}$?
Here $X$ is a closed 4-manifold, $C$ a contractible 4-manifold, $e$ an embedding of $C$ in $X$ and $\tau$ a self-diffeomorphism on $\partial C$.
\item Can we find such a triple for a family of exotic 4-manifolds with boundary?
\end{itemize}
In \cite{T1}, the author constructed a pair $(P,\varphi)$ of a non-contractible 4-manifold $P$ and a self-diffeomorphism $\varphi$ on $\partial P$, where it was called a {\it generalized cork} in the article.
In place of it, we will define a $\bar{c}h$-twist later.
The author \cite{T} and Auckly, Kim, Melvin, and Ruberman \cite{AKMR} constructed $n$-corks for any natural number $n$, where an $n$-cork is defined in Section~\ref{copl}. 
The authors in \cite{AKMR} defined a $G$-cork $(C,G)$ and gave examples.

Right after works \cite{T} and \cite{AKMR}, 
Gompf in \cite{G} constructed a new type of corks.
Let $V$ denote a local elliptic fibration over $D^2$ union three $-1$-framed 2-handles corresponding to vanishing cycles along a fiber torus:
One is a generator $\ell$ of the fiber torus; the other two are the parallel cycles orthogonal to $\ell$.
Let $X_{K}$ denote the knot-surgery of $X$ along $K$ in \cite{FS}.
Here we state Gompf's result, a cork existence theorem for a family of 4-manifolds.
\begin{thm}[Gompf \cite{G}]
\label{Gompf}
Suppose that $K_n$ is the $2n$-twist knot.
Let $X$ be a 4-manifold containing $V$.
Then there exists a ${\mathbb Z}$-effective embedding $C\hookrightarrow X$ such that $X_{K_n}=X(C,f^n)$.
\end{thm}	
A `${\mathbb Z}$-effective embedding' will be explained in Section~\ref{copl}.
The author in \cite{T2} generalized the construction of Gompf's ${\mathbb Z}$-corks to ${\mathbb Z}^n$-corks for any natural number $n$.
Masuda in \cite{Masuda} constructed $G$-corks for infinite non-abelian group $G$.
In \cite{MS}, Melvin and Schwartz gave an affirmative answer for the first question for any finite family of simply-connected, closed, 
exotic 4-manifolds.

Does a cork existence theorem for a family of non-simply-connected exotic 4-manifolds hold?
This would be a question that must be considered in the future.
As to the second question above, there are no general solutions so far.
See also Question~\ref{relvers} in Section~\ref{relresult}.

Our motivation is to find necessary conditions for infinite exotic families for the cork existence theorem to hold.
The following is our main result.
\begin{thm}
\label{main}
Let $\mathcal{X}$ be a family of closed exotic 4-manifolds with $b_2^+>1$.
If it is generated by a $G$-cork $(C,G)$, then $\mathcal{OS}_{\mathcal{X}}$ is a finite set.
\end{thm}
The finiteness of OS-invariants is the necessary condition we desire.
Furthermore, for a rational homology sphere $Y$, we will estimate the number of exotic structures which the cork twist of a cork bounding $Y$ produces using  on OS-invariants (Theorem~\ref{cor1}).
According to the condition in Theorem~\ref{main}, we obtain the following {\it cork non-existence theorem}.
This is a negative solution for the first question above.
\begin{thm}
There exists an infinite family of simply-connected, exotic closed 4-manifolds not generated by any infinite order cork.
\end{thm}
As an application we obtain a cork non-existence theorem (Theorem~\ref{notcorecork}) for families of 4-manifolds with boundary.

\subsection{Corks and plugs}
\label{copl}
Here we give general definitions of corks and plugs.
Let $Y$ be a manifold with boundary and $t$ a diffeomorphism $\partial Y\to \partial Y$.
We call a pair $(Y,t)$ a {\it twist} and $t$ a {\it twist map}.
If we have an embedding $e:Y\hookrightarrow X$ as a codimension 0 submanifold,
then a twist $(Y,t)$ induces a surgery $X\leadsto X(Y,e,t)$.

If a twist map $t$ extends to a diffeomorphism on $Y$, then we call a twist $(Y,t)$ a {\it trivial twist}, otherwise a {\it nontrivial twist}.
If $Y$ is a contractible (or non-contractible) 4-manifold, then we call the twist $(Y,t)$ a {\it $c$-twist} 
(resp. {\it $\bar{c}$-twist}).
If a twist map $t$ extends to a homeomorphism on $Y$, we call the twist $t$ an {\it $h$-twist},
otherwise an {\it $\bar{h}$-twist}.
According to Freedman's result, any $c$-twist is an $h$-twist.
The contraposition implies that any $\bar{h}$-twist is a $\bar{c}$-twist.
Hence, we have
\begin{center}
$c$-twist$\Rightarrow$ $h$-twist, $\bar{h}$-twist $\Rightarrow \bar{c}$-twist.
\end{center}
Any $c$-twist is not always the trivial twist.
It leads to the existence of corks.
If a twist $(C,g)$ is a $\bar{c}$-twist and $h$-twist, then we call it a {\it $\bar{c}h$-twist}.
An example in \cite{T1} is a $\bar{c}h$-twist.

Here we define a cork and a plug.
\begin{defn}[Cork]
\label{dfncork}
We call a nontrivial $c$-twist a {\it cork}.
For a cork $(C,f)$, if $C$ is a Stein manifold, 
then $(C,f)$ is called a {\it Stein cork}.
If $C$ is a submanifold of $X$, then the deformation $X\leadsto X(C,f)$ is called a {\it cork twist}.
\end{defn}
The pair $(C,\tau)$ in Fact~\ref{Ma} is an example of a cork.
Notice that the condition $f^2=\text{id}$ is not imposed on Definition~\ref{dfncork}.
We call $\min\{n>0|(C,\tau^n)\text{ is a trivial twist}\}$
the {\it order of a cork} $(C,\tau)$, if the minimum exists.
If the order of a cork $(C,\tau)$ is $n$, then $(C,\tau)$ is called an {\it $n$-cork}.
If the minimum does not exist, then the order of the cork is called {\it infinite}.
Such a cork is called an {\it $\infty$-cork}.
$(C,f)$ in Theorem~\ref{Gompf} is an $\infty$-cork.
Notice that the order of a cork is, in general, different from the order as a
self-diffeomorphism on the boundary.

Let $\text{Diff}^+(M)$ be the group of orientation-preserving diffeomorphisms on a manifold $M$.
We define a $G$-cork.
\begin{defn}[$G$-cork]
Let $C$ be a contractible 4-manifold with boundary and $G$ a subgroup of $\text{Diff}^+(\partial C)$.
If for any nontrivial element $g\in G$, $(C,g)$ is a cork, then $(C,G)$ is called a {\it $G$-cork}.

Let $G$ be a subgroup of $\text{Diff}^+(\partial Y)$.
If for any nontrivial element $g\in G$, $(Y,g)$ is a $\bar{c}h$-twist, then $(Y,G)$ is called a {\it $G$-$\bar{c}h$-twist}.
\end{defn}
From this definition, an $n$-cork $(C,\tau)$ with $\tau^n=\text{id}$ gives a ${\mathbb Z}_n$-cork $(C,\langle \tau\rangle)$ and for an $\infty$-cork $(C,g)$, the pair $(C,\langle g\rangle)$ is a ${\mathbb Z}$-cork.
Here ${\mathbb Z}_n$ stands for the quotient group ${\mathbb Z}/n{\mathbb Z}$.

Let $(Y,g)$ be a twist.
For an embedding $e:Y\hookrightarrow X$, if $X(Y,e,g)$ is exotic to $X$, then the embedding $e$ is called a {\it $g$-effective embedding}.

Let $G\subset \text{Diff}^+(\partial Y)$ be a subgroup.
If for any elements $g_1,g_2\in G$ with $g_1\neq g_2$, $X(Y,e,g_1)$ and $X(Y,e,g_2)$ are exotic, then $e$ is called a {\it $G$-effective embedding}.
The examples in \cite{T} or \cite{AKMR} give ${\mathbb Z}_n$-corks, $G$-corks and $G$-effective embeddings for a subgroup $G$ in $SO(4)$.

\begin{figure}[htb]
\begin{overpic}
[width=.7\textwidth]
{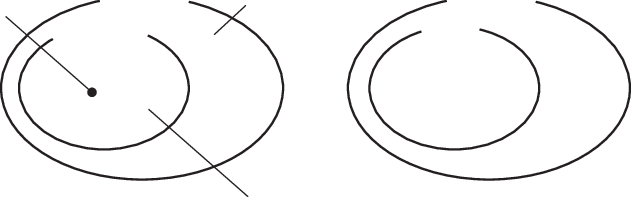}
\put(-10,30){trivial twist}
\put(10,25){$c$-twist}
\put(17,30){$h$-twist}
\put(40,-2){cork}
\put(40,30){$\bar{c}h$-twist}
\put(71,30){$\bar{h}$-twist}
\put(68,25){plug}
\end{overpic}
\caption{Relationship of various twists}
\label{ccc}
\end{figure}

We define another type of $\bar{c}$-twists.
\begin{defn}[Plug]
\label{plug}
If a twist $(P,g)$ is an $\bar{h}$-twist and has a $g$-effective embedding $P\hookrightarrow X$, then we call this twist $(P,g)$ a {\it plug}.
\end{defn}
In {\sc Figure}~\ref{ccc} we present various kinds of twists.

The result \cite{AR} shows that if $(W,F)$ is any nontrivial $h$-twist, then there is
a 4-manifold containing $W$ such that the twist gives exotic 4-manifolds.
It would be natural to consider whether any $\bar{h}$-twist has a $G$-effective embedding in a 4-manifold.
We raise the following simple question:
\begin{que}
When is an $\bar{h}$-twist $(P,g)$ a plug?
\end{que}


\subsection{Galaxy}
Let $X$ be a 4-manifold.
We define the set of smooth structures on $X$ to be gal$(X)$.
We call it the {\it galaxy} of $X$.
Let $S\subset$gal$(X)$ be a subset with $\#S>1$.
If there exist a $G$-cork $(C,G)$ and a $G$-effective embedding $e:C\hookrightarrow X$ such that $S\subset\{X(C,e,g)|g\in G\}$, then
we say that $S$ {\it is generated by a $G$-cork $(C,G)$}.
A $G$-effective embedding means that it gives an injection $G\to \text{gal}(X)$.
If $S=\{X(C,e,g)|g\in G\}$ is satisfied for a $G$-effective embedding $e$, then we call $S$ a {\it $G$-constellation}.

Fact~\ref{Ma} implies any pair of 2 points in gal$(X)$ is a $\bZ_2$-constellation
for any closed simply-connected 4-manifold $X$.
Melvin and Schwartz in \cite{MS} proved that any finite subset of gal$(X)$ is a $\bZ_n$-constellation.

Our concrete motivation of this article is the following natural question:
\begin{que}
\label{mainques}
Let $X$ be any closed oriented 4-manifold.
Is any subset $S\subset$gal$(X)$ a $G$-constellation?
\end{que}
\subsection{Finiteness condition of OS-polynomials and its examples}
To answer  Question~\ref{mainques} negatively, we use Ozsv\'ath-Sz\'abo's 4-manifold invariant $\mathcal{OS}_X$ (called an {\it OS-polynomial}).
We will explain this invariant in Section~\ref{os4}.
We set $\{\mathcal{OS}_{Y}|Y\in \mathcal{X}\}$ by $\mathcal{OS}_{\mathcal{X}}$.
Here we restate our main theorem.\\[2mm]
{\bf Theorem~\ref{main}.}
{\it Let $\mathcal{X}$ be a family of closed exotic 4-manifolds with $b_2^+>1$.
If it is generated by a $G$-cork $(C,G)$, then $\mathcal{OS}_{\mathcal{X}}$ is a finite set.
}\\

We notice that even if $\mathcal{X}\subset \text{gal}(X)$ is an infinite set, 
then $\mathcal{OS}_{\mathcal{X}}$ must be a finite set if it is generated by a cork.
As a corollary of Theorem~\ref{main}, we obtain the following.
\begin{cor}
\label{maincor}
There exist a closed oriented 4-manifold $X$ and a set $\mathcal{X}\subset \text{gal}(X)$ such that it is not generated by a $G$-cork for any group $G$.
\end{cor}

It is easy to find examples satisfying this condition.
For example, the following example illustrates Corollary~\ref{maincor}.
\begin{exm}
Let $X$ be the elliptic fibration $E(2)$ i.e., the K3 surface.
$K_n$ is the $2n$-twist knot, and $T_n$ is the $(2,2n+1)$-torus knot.
Two families $\mathcal{X}_1=\{X_{K_n}|n\in{\mathbb N}\}$ and $\mathcal{X}_2=\{X_{T_n}|n\in{\mathbb N}\}$ are two subsets of $\text{gal}(E(2))$.
Comparing the sets $\mathcal{OS}_{\mathcal{X}_1}$ and $\mathcal{OS}_{\mathcal{X}_2}$, we have the following.

$\mathcal{OS}_{\mathcal{X}_1}$ is a finite set:
$$\{1,t-1+t^{-1}\}.$$
$\mathcal{X}_1$ is a ${\mathbb Z}$-constellation by Theorem~\ref{Gompf}.

On the other hand, $\mathcal{OS}_{\mathcal{X}_2}$ is the following infinite set
\begin{equation}
\label{inf}
\left\{\sum_{k=-n}^{n}(-t)^{k}\in {\mathbb Z}_2[t^{\pm1}]\Big|n\in{\mathbb N}\right\}.
\end{equation}
Theorem~\ref{main} shows that $\mathcal{X}_2$ is not a $G$-constellation and not even generated by any $G$-cork.
\begin{figure}[htb]
\begin{overpic}
[width=.7\textwidth]
{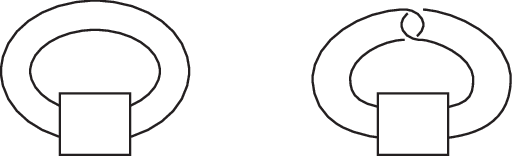}
\put(13.75,5){$n+\frac{1}{2}$}
\put(79.5,5){$n$}
\end{overpic}
\caption{The torus knot $T_{n}$ and the twist knot $K_n$}
\label{c}
\end{figure}

\end{exm}
In the same way, for $m\ge 3$, the OS-invariant of $E(m)_{T_n}$ is 
\begin{eqnarray*}
OS_{E(m)_{T_n}}&=&(t-t^{-1})^{m-2}(t^n-t^{n-1}+\cdots-t^{-n+1}+t^{-n})\bmod 2\\
&=&(t-t^{-1})^{m-3}(t^{n+1}-t^n+t^{-n}-t^{-n-1})\bmod 2.
\end{eqnarray*}
The family $\{E(m)_{T_n}|n\in {\mathbb Z}\}$ is not also generated by any $G$-cork.

This proof means that for any family $\mathcal{K}$ of knots, if $\{E(2)_{K}|K\in \mathcal{K}\}$ is generated by an infinite order cork, then $\#\{\Delta_{K}(t)\bmod 2|K\in \mathcal{K}\}$ is finite.

Similar non-existence results with respect to corks are also obtained by Yasui in \cite{Yasui}.
Our current result inspired his result.
The family of the torus knot-surgeries above is the first example not being a $G$-constellation.

\subsection{A cork for relatively exotic 4-manifolds}
\label{relresult}
Let $X$ be a 4-manifold with a {\it marking} $f:Y\to \partial X$ (diffeomorphism), 
where $Y$ is a closed 3-manifold.
We call $(X,Y,f)$ a {\it relative manifold}.
The triples $(X_1,Y,f_1)$ and $(X_2,Y,f_2)$ are {\it homeomorphic} (or {\it diffeomorphic} respectively) 
if there exists a homeomorphism (or a diffeomorphism respectively) $F:X_1\to X_2$
such that $F|_{\partial X_1}=f_2\circ  f_1^{-1}:\partial X_1\to \partial X_2$.
Consider a {\it diffeomorphism type of a relative manifold}.
Then we denote the homeomorphism (or diffeomorphism) by $(X_1,Y,f_1)\simeq (X_2,Y,f_2)$ 
(or by $(X_1,Y,f_1)\cong (X_2,Y,f_2)$ respectively).
If the equivalence class on relative 4-manifolds by homeomorphisms (or diffeomorphisms) is
called the {\it homeomorphism type} (or {\it diffeomorphism type}).

If two relative manifolds are homeomorphic but non-diffeomorphic, 
then we say that the relative manifolds are {\it exotic} or they are relatively exotic.
$(C,\tau)$ is a cork if and only if the relative 4-manifolds $(C,\partial C,\text{id})$ and $(C,\partial C,\tau)$ are 
exotic.
A twist $(X,g)$ is a nontrivial $h$-twist, if and only if $(X,\partial X,\text{id})$ and 
$(X,\partial X,g)$ are relatively exotic 4-manifolds.
A plug $(P,g)$ gives a pair of non-homeomorphic relative 4-manifolds $(P,\partial P,\text{id})$ and $(P,\partial P,g)$
which there is a $G$-effective embedding $P\hookrightarrow X$.

Here we give an open question corresponding to a relative version of Fact~\ref{Ma}.
\begin{que}
\label{relvers}
Let $(X_1,Y,f_1)$ and $(X_2,Y,f_2)$ be two relatively exotic 4-manifolds.
Then does there exist a cork $(C,\tau)$ (possibly not a 2-cork) such that $(X_2,Y,f_2)\cong (X_1(C,\tau),Y,f_1)$?
\end{que}

Next, consider the same question for any infinite family of relative 4-manifolds.
Let $(X,Y,f)$ be a relative 4-manifold.
If there exists a pair of a cork $(C,G)$ and an embedding $e:C\hookrightarrow X$ such that 
for any elements $g_1,g_2\in G$ with $g_1\neq g_2$, $(X(C,e,g_1),Y,f)$ and $(X(C,e,g_2),Y,f)$
are exotic, then we call the embedding $e$ a {\it $G$-effective embedding} in $(X,Y,f)$.

We can define a galaxy gal$(X,\partial X,f)$ for a relative 4-manifold as
$$\text{gal}(X,\partial X,f):=\{(X',\partial X,f')|(X',\partial X,f')\text{ is exotic to }(X,\partial X,f)\}/\text{diffeo}.$$
For a family $\mathcal{X}\subset$gal$(X,\partial X,f)$, 
if there exists a pair of a $G$-cork $(C,G)$ and a $G$-effective embedding $e:C\hookrightarrow X$ such that 
$\mathcal{X}\subset \{(X(C,e,g),\partial X,\text{id})|g\in G\}$, then 
we also say that $\mathcal{X}$ is a {\it $G$-constellation}.
We also say that an exotic family $\mathcal{X}$ of relative 4-manifolds {\it is generated by a $G$-cork $(C,G)$},
if the family $\mathcal{X}$ is a subset of a $G$-constellation
containing a $(X,\partial X,f)\in \mathcal{X}$ for some embedding $C\hookrightarrow X$.

In this article, we can prove a cork non-existence theorem on relative version of Corollary~\ref{maincor}.

\begin{thm}
\label{notcorecork}
There exist a relative 4-manifold $(X,\partial X, f)$ and a family $\mathcal{X}\subset \text{gal}(X,\partial X,f)$
such that $\mathcal{X}$ is not any $G$-constellation.
\end{thm}
This is proven on the basis of the idea of Theorem~\ref{main}.
Actually, the family is not even generated by any $G$-cork.

Further, in this article, we also propose ways to get corks between exotic 4-manifolds (Remark~\ref{idea}).

\subsection{OS-polynomial variations of a constellation given by a cork}
A cork should reflect certain properties of $G$-constellation.
The proof of our main theorem (Theorem~\ref{main}) gives an OS-polynomial estimate for a family to have a $G$-constellation.
In the case of $M=T[m]^{\oplus n}_{(d)}$, 
then we have 
$$|Aut_U(M)|=|SL(n,{\mathbb F})|=\prod_{k=0}^{n-1}(2^n-2^k)=2^{n(n-1)/2}\prod_{k=1}^{n}(2^{k}-1).$$
We put $N_n=|SL(n,{\mathbb F})|$.
Hence, for example, $N_1=1, N_2=6, N_3=168$, and $N_4=20160$.

\begin{thm}
\label{cor1}
Let $X$ be a closed 4-manifold with $b_2^+>1$.
For an exotic family $\mathcal{X}\subset$gal$(X)$ if there exists some group $G$ such that $\mathcal{X}$ is generated by a $G$-cork $(C,G)$ and $HF_{\text{red}}^-(\partial C)$ is a component-wise disjoint graded module of direct sum of finite dimensional module $T[p_i]^{\oplus n_i}_{(d_i)}$ $(i=1,\cdots,m)$ where $d_1,\cdots, d_m$ are distinct integers,
then we have $\#\mathcal{OS}_{\mathcal{X}}\le N_{n_1}\cdots N_{n_m}$.
\end{thm}
The definition of `component-wise disjoint' will be done in Section~\ref{actiondefinition}.
If $HF_{\text{red}}^-(\partial C)$ is component-wise disjoint, graded module and $n_1=\cdots =n_m=1$, then $\#\mathcal{OS}_{\mathcal{X}}=1$ holds.
Namely, the OS-invariants are not changed by the cork twist.
In the case of $(m,n_1)=(1,2)$, $\#\mathcal{OS}_{\mathcal{X}}$ is at most 6.

Here we observe several examples.
Other examples will be described in the later section.
\begin{exm}
In the case of the well-known Akbulut cork as in \cite[Figure 1]{AkDu}, Akbulut and Durusoy computed the Heegaard Floer homology of the boundary $\Sigma$ of the cork as
$$HF_{\text{red}}^-(\Sigma)\cong T[1]^{\oplus 2}_{(-2)}.$$
Hence, the variation of OS-invariants is at most 6.

In \cite{AkKar} Akbulut and Karakurt computed Heegaard Floer homologies of homology spheres $\partial W_n(2n+1)$ that $W_n(2n+1)$ is a contractible 4-manifold with Mazur type.
The definition is as in \cite[{\sc Figure}~3]{AkKar}.
For example since $HF_{\text{red}}^-(\partial W_1(3))\cong T[1]^{\oplus 2}_{(-2)}$ holds, we also have if a cork $C$ with $\partial C=\partial W_1(3)$ generates a family $\mathcal{X}$ of 4-manifolds, then $\#\mathcal{OS}_{\mathcal{X}}\le 6$.
Similarly, 
$$HF_{\text{red}}^-(\partial W_2(5))=T[1]_{(-2)}^{\oplus 4}\oplus T[1]_{(-4)}^{\oplus 2}\oplus T[1]_{(-12)}^{\oplus 2}.$$
Then if $(C,G)$ is a pair of a contractible 4-manifold $C$ with $\partial C=\partial W_2(5)$
and $G\subset \text{Diff}^+(\partial W_2(5))$ and a family $\mathcal{X}\subset gal(X)$ is generated by a $G$-cork $(C,G)$ then we have
$$\#\mathcal{OS}_{\mathcal{X}}\le N_4N_2^2=725760.$$
\end{exm}
It does not seem that these estimates are so sharp as to give obstructions 
enough to study corks.
We need consider the following questions for our future study.
\begin{que}
Let $Y$ be the boundary of a $G$-cork $(C,G)$.
\begin{enumerate}
\item Which diffeomorphisms on $Y$ can be a realized as an element in $G$?
\item Which automorphisms on $HF^-(Y)$ can be realized as $g_\ast$ of $g\in G$?
\end{enumerate}
\end{que}
\section{Preliminaries}
\label{Preliminaries}

\subsection{4-dimensional SW-invariant and OS-invariant.}
\label{os4}
Seiberg and Witten defined a topological invariant for a smooth 4-manifold $X$ with $b_2^+(X)>1$:
$$SW_X:\text{Spin}^c(X)\to {\mathbb Z}.$$
In the case where there are no 2-torsions in $H_1(X,{\mathbb Z})$,
the set of spin$^c$ structures is identified with $H^2(X,{\mathbb Z})$.
This gives a polynomial, {\it SW-polynomial}: 
$$\mathcal{SW}_X:=\sum_{a\in H^2(X,{\mathbb Z})}SW(a)e^a\in {\mathbb Z}[H^2(X,{\mathbb Z})].$$
Taking composition $SW_X$ with the reduction ${\mathbb Z}\to {\mathbb Z}_2$, 
we obtain a {\it mod 2 SW-polynomial}, $\mathcal{SW}^{mod 2}_X$.

Ozsv\'ath and Szab\'o defined a smooth 4-manifold invariant
$$OS_{X}:\text{Spin}^c(X)\to {\mathbb Z}_2$$
for any closed oriented 4-manifold with $b_2^+>1$.
We call the invariant $OS_{X}$ the {\it OS-invariant} and the polynomial 
$$\mathcal{OS}_X:=\sum_{a\in H^2(X,{\mathbb Z})}OS_X(a)e^a\in {\mathbb Z}_2[H^2(X,{\mathbb Z})]$$
the {\it OS-polynomial}.
It is conjectured that $\mathcal{OS}_X$ is equivalent to $\mathcal{SW}_X^{mod2}$.


We give a brief review of the OS-invariant in \cite{OS1}.
Let $(X,\frak{s})$ be a closed spin$^c$ 4-manifold with $b^+_2>1$.
There exists a 3-dimensional submanifold $N\subset W$ such that we obtain a decomposition $X=X_0\cup_NX_1$ that $b^+_2(X_i)>0$ and the induced map $H^1(N)\to H^2(X)$ is trivial.
Such a submanifold $N$ is called an {\it admissible cut}.
The spin$^c$ structure $\frak{s}$ restricts to $\frak{s}_0$ and $\frak{s}_1$ over $X_0$ and $X_1$.
We denote the restriction of $\frak{s}$ to $N$ by $\frak{t}$.
Let $X_i'$ be $X_i$ with a 4-ball deleted in the interior.
Then $X_0'$ is a cobordism from $S^3$ to $N$ and $X_1$ is a cobordism from $N$ to $S^3$.
The cobordisms $X_0'$ and $X_1'$ induce the two homomorphisms
$$F^-_{X_0',\frak{s}_0}:HF^-(S^3)\to HF^-(N,\frak{t}),$$
$$F^+_{X_1',\frak{s}_1}:HF^+(N,\frak{t})\to HF^+(S^3).$$
Since the image of $F^-_{X_0',\frak{s}_0}$ factors through $HF_{\text{red}}^-(N,\frak{t})\to HF^-(N,\frak{t})$ and 
$F_{X'_1,\frak{s}_1}^+$ factors through $HF^+(N,\frak{t})\to HF_{\text{red}}^+(N,\frak{t})$.
Composing these maps along $HF_{\text{red}}^-(N,\frak{t})\cong HF^+_{\text{red}}(N,\frak{t})$, we get $F^{\text{mix}}_{X,\frak{s}}:HF^-(S^3)\to HF^+(S^3)$.
Then $OS_{X}(\frak{s})\in{\mathbb Z}_2$ is defined to be the value satisfying
$F^{\text{mix}}_{X}(\frak{s)}(U^d\cdot\Theta^-)=OS_{X}(\frak{s})\cdot\Theta^+$, where $d=(c^2_1(\frak{s})-2\chi(X)-3\sigma(X))/4$ and $\Theta^{-}$ (or $\Theta^+$) is the maximal (or minimal) grading generator in $HF^{-}(S^3)$ (or $HF^+(S^3)$).
\subsection{Knot-surgery and SW- and OS-invariant formulas}
Let $K$ be a knot in $S^3$.
Let $X$ be a 4-manifold with a square zero embedded torus $T$.
Then the surgery 
$$X_K=[X-\nu(T)]\cup [(S^3-\nu(K))\times S^1]$$
is called a {\it (Fintushel-Stern's) knot-surgery} with respect to $K$.
The notation $\nu(\cdot)$ stands for a tubular neighborhood of a submanifold.
The gluing map is indicated in \cite{FS}.

Recall the SW-polynomial formula for the knot-surgery:
\begin{equation}
\label{SWknot}
\mathcal{SW}_{X_K}=\mathcal{SW}_{X}\cdot \Delta_K(t).
\end{equation}
Here $\Delta_K(t)$ is the Alexander polynomial of $K$.
This formula (\ref{SWknot}) means that if $\mathcal{SW}_X\neq 0$ and
 Alexander polynomials of knots $K_1,\cdots, K_n$ are mutually different, then 
the knot-surgeries $X_{K_i}$ of $X$ are mutually non-diffeomorphic.
For example, the SW-polynomial of the K3 surface, $E(2)$ is $1$.

In \cite{M}, Mark proved the knot-surgery formula of the OS-polynomial:
\begin{equation}
\label{mark}
\mathcal{OS}_{X_K}=\mathcal{OS}_{X}\cdot\Delta_K(t).
\end{equation}
This is the OS-polynomial counterpart of the formula (\ref{SWknot}).
Here we notice Sunukjian's result in \cite{Su} that the Alexander polynomial distinguishes smooth structures of Fintushel-Stern's knot-surgeries.

\subsection{$G$-action on $HF^-$}
\label{actiondefinition}
Let ${\mathbb F}$ denote the field with characteristic 2.
Let $M$ be a ${\mathbb F}[U]$-module with a ${\mathbb Q}$-grading.
The $U$-action decreases the grading by 2.

Let $\mathcal{D}(M)\subset {\mathbb Q}$ be the set of gradings of non-zero homogenous elements in $M$.
We suppose that $\mathcal{D}(M)$ is discrete and has an upper-bound.
We call a linear isomorphism $\varphi:M\to M$ a {\it $U$-automorphism} if $U\cdot \varphi(x)=\varphi(U\cdot x)$ is satisfied for any $x\in M$.
Let $Aut_U(M)$ be the group of $U$-automorphisms over $M$ with grading-preserving.

Here $d^-(Y)$ is the maximal grading of the image of $HF^-(Y)\to HF^\infty(Y)$ and it is computed by $d^-(Y)=d(Y)-2$, where $d(Y)$ is called the {\it correction term}.
If $Y$ has a contractible bound, then $d^-(Y)=-2$ holds.

Let $T^-_{(d)}$ denote ${\mathbb F}[U]$ with the top (absolute) grading $d\in {\mathbb Q}$.
We call it a {\it free ${\mathbb F}[U]$-module}.
Let $T[m]_{(d)}$ denote $T_{(d)}^-/(U^m)$ with the top grading $d$.
This module is a finite dimensional ${\mathbb F}$-vector space.
Simply, we say that such a module is {\t finite dimensional}.
Here we define `component-wise disjoint'.
\begin{defn}
Let $M_i$ ($i=1,2,\cdots$) be a ${\mathbb F}[U]$-module isomorphic to $T[p_i]_{(d_i)}^{\oplus n_i}$, and $(p_i,d_i)\neq (p_j,d_j)$ for $i\neq j$.
For such a $M_i$, we put $M=\oplus_{i}M_i$.
Then $M_i$ ($i=1,2,\cdots$) is called a {\it component} of $M$.
If for $i\neq j$, $\mathcal{D}(M_i)\cap \mathcal{D}(M_j)=\emptyset$ is satisfied, then we say that $M_i$ and $M_j$ are {\it disjoint}.
For any $i\neq j$, if $M_i$ and $M_j$ are disjoint, then we call $M=\oplus_i M_i$ {\it component-wise disjoint}.
\end{defn}
It follows that we have
\begin{lem}
\label{keylem}
Let $M$ be a ${\mathbb Q}$-graded, finite dimensional ${\mathbb F}[U]$-module. 
Then $M$ is isomorphic to $\oplus_{i=1}^mT[n_i]_{(d_i)}$.
\end{lem}
If a non-zero homogenous element $x\in M$ satisfies $x\neq U\cdot y$ for any $y\in M$, then we call $x$ a {\it maximal element}.
Notice that if a graded module $M$ is finite dimensional, then $\mathcal{D}(M)$
is bounded.

\begin{proof}
We take a basis $\{x_1,\cdots,x_m\}$ of $\ker(U)\subset M$ each of which has a homogenous grading.
Let $y_i$ be a maximal element satisfying $x_i=U^{p_i-1}\cdot y_i$ for any $i$.
Here we may assume that $y_i$ is also homogenous and that $p_1\ge p_2\ge \cdots\ge p_m$.
We assume that the grading of $y_i$ is $d_i$.
We show $\mathcal{B}:=\{U^{j}\cdot y_i|0\le j< p_i, 1\le i\le m\}$ is a basis of $M$ as a ${\mathbb F}$-vector space.


We prove linearly independence.
For $c_{i,j}\in {\mathbb F}$
$$\sum_{i=1}^m\sum_{j=0}^{p_i-1}c_{i,p_i-j}U^j\cdot y_i=0,$$
then we have
\begin{eqnarray*}
0&=&U^{p_1-1}\left(\sum_{i=1}^m\sum_{j=0}^{p_i-1}c_{i,p_i-j}U^j\cdot y_i\right)\\
&=&\sum_{p_1=p_i}c_{i,p_i}U^{p_i-1}\cdot y_i=\sum_{p_1=p_i}c_{i,p_i}x_i
\end{eqnarray*}
From the linear independence of $\{x_i\}$, we have $c_{i,p_i}=0$ for $i$ with $p_1=p_i$.
In the case of $p_1=p_2=\cdots =p_{n_1}>p_{n_1+1}$, 	
multiplied by $U^{p_1-2}, U^{p_1-3},\cdots, U^{p_{n_1+1}}$ in the same way,
$c_{i,j}=0$ holds for $1\le i\le n_1$ and $p_{n_1+1}< j\le  p_1$.
Here $n_1+1$ is the maximum $i$ satisfying $p_i<p_1$. 
Next, iterating the same thing for $p_{n_1+1}=p_{n_1+2}=\cdots p_{n_1+n_2}>p_{n_1+n_2+1}$,
we have $c_{i,j}=0$ holds for $1\le i\le n_1+n_2$ and $p_{n_2+1}< j\le p_{n_1}$.
Here $n_1+n_2+1$ is the maximum $i$ satisfying $p_i<p_{n_1+1}$.
By induction, we can finish this process until $j=1$, hence, we have $c_{i,j}=0$ for any $(i,j)$.
Hence, $\mathcal{B}$ is linearly independent.

For any non-zero $y\in M$, there exists non-negative integer $k$ such that $U^ky\neq 0$ and $U^{k+1}y=0$.
If $k=0$, then clearly $y$ is a linear combination of $\mathcal{B}$.
We suppose that $k>0$.
Then $U^ky=a_1x_1+\cdots+a_mx_m=a_1U^{p_1-1}y_1+\cdots+a_mU^{p_m-1}y_m$
for $p_m>1$.
Thus, $U(U^{k-1}y-(a_1U^{p_1-2}y_1+\cdots+a_mU^{p_m-2}y_m))=0$.
Thus $U^{k-1}y-(a_1U^{p_1-2}y_1+\cdots+a_mU^{p_m-2}y_m)=b_1x_1+b_2x_2+\cdots+b_mx_m$ for some $b_i\in {\mathbb F}$.
This means that $U^{k-1}y$ is a linear combination of $\mathcal{B}$.
If $k=1$, then $y$ is a linear combination of $\mathcal{B}$.
Otherwise we can see that $U^{k-2}y$ is a linear combination of $\mathcal{B}$ in the same way.
Inductively, we see that $y$ is a linearly combination of $\mathcal{B}$ in the end.
As a result, $\mathcal{B}$ is a basis of $M$.
This means that $M$ is isomorphic to $\oplus_{i=1}^mT[p_i]_{(d_i)}$.
\qed
%
%
\end{proof}
For a homology sphere $Y$, $HF^-(Y)$ is isomorphic to the direct sum 
$HF^-(Y)=T^-_{(d^-(Y))}\oplus HF_{\text{red}}^-(Y)$, where $HF_{\text{red}}^-(Y)$ 
is the kernel of the natural homomorphism $HF^-(Y)\to HF^\infty(Y)$ and is a finite dimensional ${\mathbb F}[U]$-module.

\begin{lem}
\label{keyfinite}
There is an isomorphism $Aut_U(HF^-(Y))\cong Aut_U(HF^-_{\text{red}}(Y))$.
In particular $Aut_U(HF^-(Y))$ is a finite group.
\end{lem}
\begin{proof}
An automorphism $\varphi \in Aut_U(HF^-(Y))$ induces the following commutative diagram:
$$
\xymatrix{
 HF^-(Y) \ar[d]^{\varphi}\ar[r] & HF^\infty(Y) \ar[d]^{\varphi_\infty=\text{id}} \\
 HF^-(Y) \ar[r] & HF^\infty(Y)
}
$$
where $\varphi_\infty\in Aut_U(HF^\infty(Y))$ is induced.
Since $\varphi$ preserves the grading, $\varphi_\infty=\text{id}$ holds.
This means that $\varphi$ is the identity on $T^-_{(d^-(Y))}$.
Thus, the group isomorphism $Aut_U(HF^-(Y))\cong Aut_U(HF_{\text{red}}^-(Y))$ is satisfied.

Since $HF_{\text{red}}^-(Y)$ is a finite dimensional ${\mathbb F}[U]$-module then 
$Aut_U(HF^-(Y))$ is a finite group.
In fact, this group is a subgroup of $SL(N,{\mathbb F})$ for a sufficiently large integer $N$.
\qed
\end{proof}
For example, in the case of $M=T[m]_{(d)}^{\oplus n}$, we have $Aut_U(M)=SL(n,{\mathbb F})$ by considering where the maximal elements $M/U\cdot M={\mathbb F}^n_{(d)}$ move.

%

We prove the following proposition before proofs of our theorems.
\begin{prop}
\label{mainp}
Let $\mathcal{X}$ be a family of exotic closed oriented 4-manifolds with $b_2^+>1$.
If $\mathcal{X}$ is generated by twists of a $G$-cork, then 
the OS-polynomials $\{\mathcal{OS}_{X}|X\in \mathcal{X}\}$ must be a finite set.
\end{prop}
\begin{proof}
Let $X$ be a closed oriented 4-manifold with $b^+_2(X)>1$ and $(C,G)$ a $G$-cork.
Let $C\hookrightarrow X$ be an $G$-effective embedding satisfying $\mathcal{X}\subset \{X(C,g)|g\in G\}$.

Thus, this induces a $U$-isomorphism on $HF^-(\partial C)$.
The $U$-isomorphism $g_\ast$ for $g\in G$ preserves the absolute grading of $HF^-(\partial C)$.
In fact, the shift of the grading can be computed by $(c^2_1(\frak{t})-3\sigma-2\chi)/4$.
Since $c^2_1(\frak{t})$, $\sigma$ and $\chi$ are zero for the cylinder $I\times \partial C$, the grading is preserved.
Here $\frak{t}$ is a spin$^c$-structure on $I\times \partial C$.
This means a homomorphism $G\to Aut_U(HF^-(\partial C))$.
Due to Lemma~\ref{keyfinite}, $Aut_U(HF^-(\partial C))$ is a finite group.

Here for any $g \in G$ we consider $X(C,g)$ as a composition of  three 4-manifolds $\{C,P,\tilde{V}\}$,
where  $X(C,g)=[C\cup_{g} P]\cup \tilde{V}$ and $P$ is a cobordism $\partial C \overset{P}{\to} \partial \tilde{V}$.
Furthermore, we assume that $N=\partial \tilde{V}$ is an admissible cut of $X(C,g)$.
The existence of this cut $N$ is guaranteed in \cite{OS}.
Deleting two 4-balls in the interior in $X(C,g)$, we give a composition $W_g$ of three cobordisms:
$$W_g:S^3\overset{C_0}{\to}\partial C\overset{P}{\to}N\overset{V}{\to}S^3 ,$$
where $V$ is $\tilde{V}$ with a 4-ball deleted.
Here the mixed invariant on $W_g$ is computed as follows:
$$F_{W_g,\frak{s}_0}^{\text{mix}}=F^+_{V,\frak{s}_3}\circ F_{P,\frak{s}_2}^-\circ g_\ast\circ F_{C_0,\frak{s}_1}^-:HF^-(S^3)\to HF^+(S^3),$$
where $\frak{s}\in \text{Spin}^c(X(C,g))$, $\frak{s}|_{W_g}=\frak{s}_0$, $\frak{s}_0|_{C_0}=\frak{s}_1$, $\frak{s}_0|_{P}=\frak{s}_2$ and $\frak{s}_0|_{V}=\frak{s}_3$.

As mentioned above, $g_\ast$ is regarded as a finite group action on $HF^-(\partial C)$, i.e., $\{g_\ast|g\in G\}$ is finite.
The mixed invariants $\{F^{\text{mix}}_{W_g,\frak{s}_0}|g\in G\}$ are also finite.
Thus the set $\{\mathcal{OS}_{X(C,g)}|g\in G\}$ of OS-invariants is also finite.
\qed
\end{proof}
Therefore, a family $\mathcal{X}$ generated by a single $G$-cork has to have finite variations of OS-invariants.

\section{Proofs of theorems~\ref{main} and ~\ref{notcorecork}}
\subsection{Proof of Main theorem~\ref{main}.}
Let $\mathcal{X}$ be a set of infinite smooth structures of closed oriented 4-manifolds with $b_2^+>1$.
If OS-invariants of $\mathcal{X}$ have infinite variations, then it is not generated by a single $G$-cork by using Proposition~\ref{mainp}.
\qed

{\bf Proof of Theorem~\ref{cor1}.}
Let $(\mathcal{C},G)$ be a $G$-cork generating $\mathcal{X}\subset gal(X)$
for some 4-manifold $X$ with $b_2^+(X)>1$.
We put $Y=\partial \mathcal{C}$.
Any element $g\in G$ induces
an action $g_\ast$ on $HF^-(Y)$.
From assumption, there exists an isomorphism $HF^-_{\text{red}}(Y)\cong \oplus_{i=1}^mM_i$, where $\{M_i=T[p_i]_{(d_i)}^{\oplus n_i}|i=1,\cdots, m\}$ are component-wise disjoint.

Then since actions $g_\ast$ are grading-preserving, by Lemma~\ref{keylem} the actions on $HF^-(Y)$ are realized as elements in 
$\prod_{i=1}^mSL(n_i,{\mathbb F})$.
Hence the cardinality of $\mathcal{OS}_{\mathcal{X}}$ is less than or equal to $N_{n_1}\cdots N_{n_m}$.
\hfill$\Box$

\subsection{Examples}
\begin{exm}
Consider $Y_{n}=\Sigma(2,3,12n+1)$.
Each of $Y_n$ has Rohlin invariant zero.
In fact, $Y_{1}$ is the boundary of a contractible 4-manifold.
Then $HF^-(Y_{n})=T^-_{(-2)}\oplus T[1]_{(-2)}^{\oplus 2n}$ and thus $HF_{\text{red}}^-(Y_{n})\cong {\mathbb F}_{(-2)}^{\oplus 2n}$.
Suppose that there exist a contractible 4-manifold $C_{n}$ with $Y_n=\partial C_n$ and a $c$-twist $\tau_n:Y_{n}\to Y_{n}$ generating a family $\mathcal{X}_n$ of closed 4-manifolds with $b_2^+>1$.
Then $\#\mathcal{OS}_{{\mathcal X}_n}\le N_{2n}$.
Heegaard Floer homology groups of $\Sigma(2,5,7)$, $\Sigma(3,4,5)$ are all isomorphic to that of $Y_1$.
Hence, in each case of them if there is a similar family $\mathcal{X}$ of 4-manifolds, then $\#\mathcal{OS}_{\mathcal X}\le N_2=6$.
\end{exm}

\begin{exm}
Rustamov in \cite{Ru} computed the Heegaard Floer homology groups of $Z_n=\Sigma(2,2n+1,4n+3)$.
Let $q_i$ be $-i(i+1)-2$ and let $p_i$ be $1,1,2,2,3,3,\cdots$ ($i=1,2,\cdots$)
Then 			
$$HF^-_{\text{red}}(Z_n)=T[p_n]_{(-2)}\oplus_{i=1}^{n-1}T[p_i]_{(q_{n-i})}^{\oplus 2}$$


Suppose $Z_n$ is the boundary of a cork and an exotic family $\mathcal{X}_n$ of closed 4-manifolds with $b_2^+>1$ which is generated by the cork.
In the case of $n=1,2$, since 
$$HF^-_{\text{red}}(Z_n)=\begin{cases}T[1]_{(-2)}&n=1\\T[1]_{(-2)}\oplus T[1]_{(-4)}^{\oplus 2}&n=2,\end{cases}$$
$\#\mathcal{OS}_{\mathcal{X}_n}\le 6^{n-1}$ holds using Theorem~\ref{cor1}.
In the case of $n=3$, $HF^-_{\text{red}}(Z_3)=T[2]_{(-2)}\oplus T[1]_{(-4)}^{\oplus 2}\oplus T[1]_{(-8)}^{\oplus 2}$.
This is not component-wise disjoint.
However, since $M_1=T[2]_{(-2)}\oplus T[1]_{(-4)}^{\oplus 2}$ and $M_2= T[1]_{(-8)}^{\oplus 2}$ are disjoint, we can obtain the following estimate 
$\mathcal{OS}_{\mathcal{X}_3}\le N_3N_2=1008$.

\end{exm}

\subsection{Non-existence theorem for relative 4-manifolds}
A $\bar{c}h$-twist gives relatively exotic 4-manifolds.
We give a family of relatively exotic 4-manifolds which is not generated by
a cork.

\begin{figure}[htbp]
\begin{center}
\includegraphics[width=.5\textwidth]{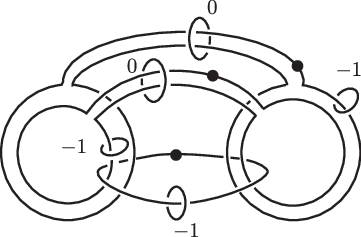}
\caption{$P$.}
\label{Pdi}
\end{center}
\end{figure}

\begin{figure}[htbp]
\begin{center}
\includegraphics[width=.8\textwidth]{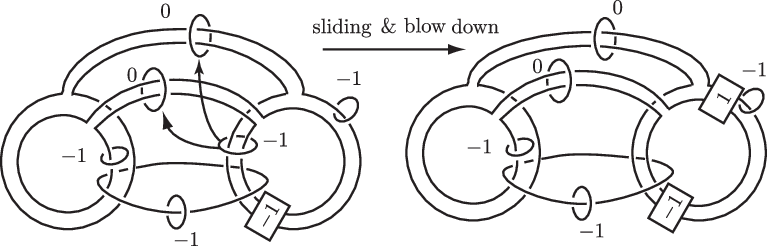}
\caption{The diffeomorphism $\varphi$.}
\label{phi}
\end{center}
\end{figure}

{\bf Proof of Main theorem~\ref{notcorecork}.}
Let $(P,\varphi)$ be the plug defined in \cite{T1}.
Namely, $P$ and $\varphi$ are described in {\sc Figure}~\ref{Pdi} and \ref{phi} respectively.

Taking $\psi=\varphi^2$, we obtain a $\bar{c}h$-twist $(P,\psi)$ by \cite{T1}.
Lemma 3.2 in \cite{T1} says that $\psi$ induces the trivial map on the homology group.

$(P,\varphi)$ changes any crossing for Fintushel-Stern's knot-surgery.
Namely, if a knot $K^+$ is one crossing change of a knot $K$,
then the plug twist changes $X_{K}$ to $X_{K^+}$.
As proven in \cite{T1}, there exists an embedding $P\hookrightarrow E(2)$ such that the twist gives the equation $E(2)(P,\psi^n)=E(2)_{T_{2n}}$.
Thus $(P,\langle \psi\rangle)$ is a ${\mathbb Z}$-$\bar{c}h$-twist.

Suppose there exists a ${\mathbb Z}$-cork $(\D,\langle f\rangle)$ generating relatively exotic 4-manifolds $\{(P,\partial P, \psi^n)|n\in {\mathbb Z}\}$.

The embedding $e:\D\hookrightarrow P$ gives the following inclusion
$$\{(P,\partial P,\psi^n)|n\in {\mathbb Z}\}\subset \{(P(\mathcal{D},f^k),\partial P,id)|k\in {\mathbb Z}\}.$$
By performing the twists for $E(2)$, we naturally obtain 
$$\{E(2)(P,\psi^n)|n\in {\mathbb Z}\}=\{E(2)_{T_{2n}}|n\in {\mathbb Z}\}\subset \{E(2)(\D,f^k)|k\in {\mathbb Z}\}$$ 
The OS-polynomials of these are
$$\{\Delta_{T_{2n}}(t)\bmod 2|n\in {\mathbb Z}\}\subset\{\mathcal{OS}_{E(2)(\D,f^k)}|k\in {\mathbb Z}\}.$$
The left hand side is an infinite set, but the right hand side is a finite set by Main theorem~\ref{main}.
This is a contradiction.
Therefore, there are no ${\mathbb Z}$-corks generating $\{(P,\partial P,\psi^n)|n\in {\mathbb Z}\}$.
\qed
\section{A plug $(Q,\eta)$ in $(P,\varphi)$.}
\subsection{A Stein plug $(Q,\eta)$}
\label{Steinplug}
\begin{figure}[htbp]
\begin{center}
\includegraphics{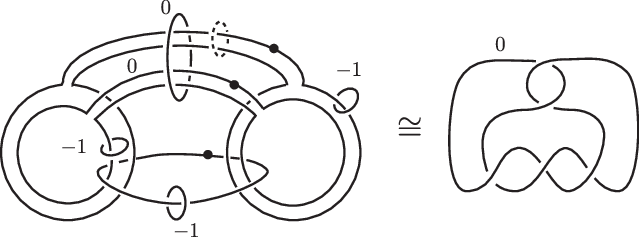}
\caption{$Q$ is the diffeomorphic to a 4-ball with a 0-framed 2-handle along $5_2$.}
\label{dif}
\end{center}
\end{figure}

In this section we define a Stein plug $(Q,\eta)$ with $b_2(Q)=1$.
Let $Q$ be the left picture in {\sc Figure}~\ref{dif} other than the dashed unknot.
By using an easy handle calculus, the $Q$ is diffeomorphic to a 4-ball with 0-framed 2-handle along $5_2$,
in particular, $\partial Q\cong S^3_0(5_2)$.
Attaching a 0-framed 2-handle along the dashed circle in the picture, we can construct $P$ in {\sc Figure}~\ref{Pdi}.
Hence $Q$ is a submanifold of $P$.
The performance of handle calculus for the map $\varphi:\partial P\to \partial P$ in \cite{T1} can induce the map $\partial Q\to \partial Q$ naturally.
Compare {\sc Figure}~\ref{phi} and {\sc Figure}~\ref{dif}.
We define the diffeomorphism to be $\eta$.
We prove the following:
\begin{prop}
\label{b1stein}
The plug $(Q,\eta)$ is a Stein plug and generates relative 4-manifolds
$\{(P,\partial P,\varphi^n)|n\in {\mathbb Z}\}$.
In particular $(Q,\eta^2)$ generates $\{(P,\partial ,\psi^n)|n\in {\mathbb Z}\}$.
\end{prop}
\begin{proof}
The handle diagram of $Q$ can be reduced to $5_2$ with framing $0$.
{\sc Figure}~\ref{reducedplug} gives the maximal Thurston-Bennequin invariant 1 Legendrian knot on $5_2$.
Thus, Gompf's result in \cite{G3} means that the presentation gives a Stein structure on $Q$.

\begin{figure}[htbp]
\begin{center}
\includegraphics[width=.3\textwidth]{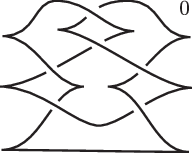}
\caption{A Stein structure on $Q$.}
\label{reducedplug}
\end{center}
\end{figure}
Attaching a 0-framed 2-handle for the diffeomorphism $\eta:\partial Q\to \partial Q$, we obtain $\varphi:\partial P\to \partial P$.
This effect of the diffeomorphism over the diagram corresponds to adding $1$-twist and $-1$-twist as in {\sc Figure}~\ref{phi}.
Since the twist via $\eta$ realizes $\varphi$, we obtain
\begin{equation}
\label{PQrel}
(P(Q,\eta^n),\partial P,\text{id})=(P,\partial P,\varphi^n)
\end{equation}

Since the twisted double $Q\cup_{\eta^k}(-Q)$ by $(Q,\eta^k)$ is homeomorphic to ${\mathbb C}P^2\#\overline{{\mathbb C}P^2}$ ($k$: odd) and
$S^2\times S^2$ ($k$: even) by easy calculation.
The diffeomorphism $\eta$ cannot extend to a homeomorphism on $Q$ by \cite{B}.
Thus $(Q,\eta^{2k+1})$ is a plug and $(Q,\eta^{2k})$ is a $\bar{c}h$-twist.

The equality~(\ref{PQrel}) means that for any unknotting number $1$ knot $K$,
there exists an embedding $Q\hookrightarrow E(2)$ such that $E(2)(Q,\eta)=E(2)_{K}$.
Thus $(Q,\eta)$ is an infinite order plug.\qed
\end{proof}
The diagram presentation in {\sc Figure}~\ref{reducedplug}
is well-known as the Chekanov-Eliashberg knot.

This proposition says that $(Q,\eta)$ is a minimal $b_2$ plug which gives rise to any crossing change for Fintushel-Stern's knot-surgery.
If there exists a plug with $b_2=0$ with the same property, then the action on the Heegaard Floer homology on the boundary, which is a rational
homology sphere, should admit infinite order.
Some iterated crossing changes give an effective action on the set of Alexander polynomials with ${\mathbb F}$-coefficients.
For example for $(2,2n+1)$-torus knots, there exists a crossing that they are produced by an iteration of the crossing change such that the iteration gives an affective action on the Alexander polynomials with ${\mathbb F}$-coefficients.

\begin{rmk}
\label{idea}
The equality~(\ref{PQrel}) implies the new method to prove that a twist is a cork as follows.
Let $(\mathcal{C},g)$ be an $\infty$-cork.
Suppose that there is an embedding $\mathcal{D}\hookrightarrow \mathcal{C}$ and a twist $h: \partial \mathcal{D}\to \partial\mathcal{D}$ such that for any integer $i$
$$(\mathcal{C}(\mathcal{D},h^n),\partial \mathcal{C},\text{id})=(\mathcal{C},\partial \mathcal{C},g^n)$$
is satisfied.
Then we say that the twist $(\mathcal{D},h)$ induces the twist $(\mathcal{C},g)$.
If a twist induces a cork, then the twist is a nontrivial.
As a result we have
\begin{observation}
To find a contractible 4-manifold which induces a cork leads to finding a new cork.
\end{observation}
This is a useful observation to find a cork contained in a known cork.
\end{rmk}

Lastly, we give a remark and question.
\begin{rmk}
It can be easily checked that $Q\cup_{\eta}(-Q)$ and $Q\cup_{\eta^2}(-Q)$ are diffeomorphic to ${\mathbb C}P^2\#\overline{{\mathbb C}P^2}$ and 
$S^2\times S^2$ respectively.
\end{rmk}
Here we give a natural question.
\begin{que}
Is any twisted double $Q\cup_{\eta^k}(-Q)$ diffeomorphic to ${\mathbb C}P^2\#\overline{{\mathbb C}P^2}$ or $S^2\times S^2$?
\end{que}

\section*{Acknowledgements}
The author is deeply grateful for some discussion with Robert Gompf, who gave me suggestions and hints to get main theorem.
The main results were inspired by the discussion.
The author thanks him for giving many helpful comments on an earlier version of the manuscript of this paper.
Furthermore, the author thanks for Kouichi Yasui giving me an opportunity to talk about this result in the Hiroshima University Topology-Geometry seminar in 2016 April.

\end{document}